\documentclass[12pt]{amsart}

\usepackage{amsfonts}
\usepackage{amssymb}
\usepackage{mathrsfs} 
\usepackage[letterpaper, left=2.5cm, right=2.5cm, top=2.5cm,
bottom=2.5cm,dvips]{geometry}
\usepackage{verbatim}
\usepackage[T1]{fontenc}
\usepackage{graphicx}
\usepackage{amsmath}
\usepackage{float}
\usepackage{polski}
\usepackage{diagbox,booktabs}
\usepackage[english]{babel}
\usepackage{pgf,tikz}
\usepackage{subcaption}
\usepackage{color}
\newtheorem{theorem}{Theorem}
\theoremstyle{plain}

\newtheorem{conjecture}{Conjecture}
\newtheorem{corollary}{Corollary}

\newtheorem{lemma}{Lemma}

\numberwithin{equation}{section}

\def\ka #1{\mathscr{#1}}
\def\kal #1 #2{\mathscr{#1}^{#2}}

\def\binom #1#2{{#1\choose #2}}

\def\skalar #1#2{\langle #1,#2\rangle}

\usepackage{diagbox} 

\usepackage{array,multirow}

\usepackage{xcolor}

\usepackage[backref]{hyperref}
\hypersetup{
	colorlinks,
	linkcolor={blue!60!black},
	citecolor={red!60!black},
	urlcolor={red!60!black}
}

\DeclareMathOperator{\dist}{dist}
\DeclareMathOperator{\diam}{diam}

\begin{document}
	\title{Neighborly boxes and strings with jokers; constructions and asymptotics}
	
		\author{Jaros\l aw Grytczuk}
\address{Faculty of Mathematics and Information Science, Warsaw University
	of Technology, 00-662 Warsaw, Poland}
\email{jaroslaw.grytczuk@pw.edu.pl}

\author{Andrzej P. Kisielewicz}
\address{Wydzia{\l} Matematyki, Informatyki i Ekonometrii, Uniwersytet Zielonog\'orski, ul. Podg\'orna 50, 65-246 Zielona G\'ora, Poland}
\email{A.Kisielewicz@wmie.uz.zgora.pl}

\author{Krzysztof Przes\l awski}
\address{Wydzia{\l} Matematyki, Informatyki i Ekonometrii, Uniwersytet Zielonog\'orski, ul. Podg\'orna 50, 65-246 Zielona G\'ora, Poland}
\email{K.Przeslawski@wmie.uz.zgora.pl}

\begin{abstract}
	We study families of axis-aligned boxes in a $d$-dimensional Euclidean space $\mathbb{R}^d$ whose placement is restricted by bounds on the dimension of their pairwise intersections. More specifically, two such boxes in $\mathbb{R}^d$ are said to be \emph{$k$-neighborly} if their intersection has dimension at least $d-k$ and at most $d-1$. The maximum number of pairwise $k$-neighborly boxes in $\mathbb{R}^d$ is denoted by $n(k,d)$. It is known that $n(k,d)=\Theta(d^k)$, for fixed $1\leqslant k\leqslant d$, however, exact formulas are known only in three cases: $k=1$, $k=d-1$, and $k=d$. In particular, the equality $n(1,d)=d+1$ is equivalent to the famous theorem of Graham and Pollak concerning partitions of complete graphs into complete bipartite graphs.
	
	In our main result we give a new construction of families of $k$-neighborly boxes which improves the lower bound for $n(k,d)$ when $k$ is close to $d$. Together with some recent upper bounds on $n(k,d)$, it gives the asymptotic equality $n(d-s,d)\thicksim\frac{2^s+1}{2^{s+1}}\cdot2^d$, for every fixed $s\leqslant d/2$. In our constructions we use a familiar interpretation of the problem in the language of Hamming cubes represented by binary strings with a special blank symbol, called \emph{joker}.
\end{abstract}
\maketitle
\begin{section}{Introduction}
	A \emph{box} in the space $\mathbb{R}^d$ is any axis-parallel $d$-dimensional cuboid. The intersection of two boxes is a again a box of possibly smaller dimension. Two boxes are said to be \emph{$k$-neighborly}, for a fixed positive integer $k\leqslant d$, if their intersection is a box of dimension strictly smaller than $d$ but not less than $d-k$. Let $n(k,d)$ denote the maximum size of a family of pairwise $k$-neighborly boxes in $\mathbb{R}^d$.
	
	The exact values of $n(k,d)$ are known only in the following three cases. In the case $k=d$ we trivially have $n(d,d)=2^d$. For $k=1$, Zaks \cite{Zaks3} proved that $n(1,d)=d+1$ by relating the problem to the celebrated theorem of Graham and Pollak \cite{GP} on bipartite partitions of cliques. Finally, for $k=d-1$, the equality $n(d-1,d)=3\cdot 2^{d-2}$ was derived only recently in \cite{AlonGKP}. The first general bounds were determined by Alon \cite{Alon} who proved that for all $1\leqslant k\leqslant d$:
	\begin{equation}
		\frac{1}{k^k}\cdot d^k\leqslant n(k,d)\leqslant \frac{2\cdot(2e)^k}{k^k}\cdot d^k.
	\end{equation}
	These bounds were recently improved in \cite{AlonGKP} and \cite{ChengWXY}. In particular, in \cite{AlonGKP}, the following lower bound for $n(k,d)$ was obtained:
	\begin{equation}
		n(k,d)\geqslant (1-o(1))\frac{d^k}{k!}.
	\end{equation}
	Moreover, the following conjecture was posed in \cite{AlonGKP}.
	\begin{conjecture}\label{Conjecture Limit}
		For every fixed integer $k\geqslant 1$, there exists a real number $\gamma_k$ such that
		\begin{equation}
			\lim_{d\rightarrow \infty}\frac{n(k,d)}{d^k}=\gamma_k.
		\end{equation}
	\end{conjecture}
	By the result of Zaks, $\gamma_1=1$, but for every $k\geqslant 2$ the conjecture is widely open. It is however tempting to guess that perhaps $\gamma_k=\frac{1}{k!}$ for every $k\geqslant 1$. Indeed, in \cite{AlonGKP} we made another supposition, which (if true) would imply this guess.
	\begin{conjecture}\label{Conjecture Pascal}
		For every $1\leqslant k\leqslant d$,
		\begin{equation}
			n(k,d)\leqslant n(k-1,d-1)+n(k,d-1),
		\end{equation}
	accepting where necessary that $n(0,d)=1$ for all $d\geqslant 1$, and $n(k,d)=2^d$ for $k>d$.
	\end{conjecture}
For instance, for $k=2$ this would give $n(2,d)\leqslant d^2/2+O(d)$, implying that $\gamma_2=1/2$. In general, by induction and the well known formula for the sum $1^k+2^k+\cdots +d^k$, one easily gets $n(k,d)\leqslant d^k/k!+O(d^{k-1})$, which shows that Conjecture \ref{Conjecture Pascal} implies Conjecture \ref{Conjecture Limit} with $\gamma_k=1/k!$.
Our main result in this paper reads as follows.
\begin{theorem}
	\label{main}
	Let $s$ and $d$ be positive integers, with $s\leqslant d/2$ and $d$ even. Then
	$$
	n(d-s,d)\geqslant   \sum_{i=0}^{\frac d 2 -s} \binom {d-s} i +\sum_{k=1}^{s-1}  2^k\binom {d-s}{\frac d 2-s +k} + 2^s\sum_{i=\frac d 2}^{d-s} \binom {d-s} i .
	$$
\end{theorem} 

The proof is by a construction conducted in the related setting of strings (see \cite{AlonGKP} or \cite{GKP}). It is given in the next section. In the last section we derive, as consequence of this result and the upper bound found in \cite{ChengWXY}, the following asymptotic estimates.
	\begin{theorem}\label{asymptotics} Let $s$ and $d$ be positive integers, with $s\leqslant d/2$. Then
	$$
	\left (\frac 1 2 +\frac 1 {2^{s+1}}\right)2^d -o(2^d) \leqslant n(d-s,d) \leqslant \left(\frac 1 2 +\frac 1 {2^{s+1}}\right)2^d.
	$$
	In consequence, for every fixed $s\leqslant d/2$, we have
	$$\lim_{d\rightarrow \infty}\frac{n(d-s,d)}{2^d}=\frac{2^s+1}{2^{s+1}}.$$
\end{theorem}

	\end{section}
	
\subsection{Binary strings with jokers}
The problem of neighborly boxes in $\mathbb{R}^d$ can be encoded in a purely combinatorial setting using strings over alphabet with just three symbols. Let $J=[0,1]$ be the unit segment and let $J_0=[0,1/2]$ and $J_1=[1/2,1]$ be its left and right half, respectively. A \emph{normalized} box is a $d$-dimensional cuboid of the form $A=A_1\times\cdots \times A_d$, where $A_i\in \{J_0,J_1,J\}$ for all $i=1,2,\ldots, d$.

Suppose that we are given two normalized boxes, $A=A_1\times\cdots \times A_d$ and $B=B_1\times\cdots \times B_d$. If for some fixed coordinate $i$, we have $\{A_i,B_i\}=\{J_0,J_1\}$, then we say that $A$ and $B$ \emph{pass} each other in dimension $i$. Otherwise, we say that $A$ and $B$ \emph{overlap} in dimension $i$. Clearly, the intersection $A\cap B$ is a cuboid whose dimension equals exactly the number of dimensions in which $A$ and $B$ overlap. For instance, if $A$ and $B$ overlap in exactly $d-1$ dimensions, or the same, if they pass each other in exactly one dimension, then $A$ and $B$ are neighborly. In general, two normalized boxes are $k$-neighborly if and only if they pass in at least one and at most $k$ dimensions.

It is not hard to imagine that any family of boxes in $\mathbb{R}^d$ can be transformed to a family of normalized boxes preserving dimensions of all intersecting pairs. Therefore in investigating the function $n(k,d)$ one may only restrict to normalized boxes.

To further simplify the setting, let $S=\{0,1,\ast\}$ be an alphabet consisting of two binary digits and one special symbol called \emph{joker}. Let $S^d$ be the set of all strings of length $d$ over $S$. Clearly, a normalized box $A=A_1\times\cdots \times A_d$ can be identified with a string $u=u_1u_2\cdots u_d$ so that $u_i=0$ if $A_i=J_0$, $u_i=1$ if $A_i=J_1$, and $u_i=\ast$ if $A_i=J$.  

		
			
%

The \emph{distance} between two strings $u,v\in S^d$ is defined as the number of positions where they differ, but none of them is occupied by a joker. It is denoted by $\dist(u,v)$. More formally, if $u=u_1\cdots u_d$ and $v=v_1\cdots v_d$, then
$$
\dist(u,v)=|\{1\leqslant i\leqslant d: u_i\neq v_i \text{ and }u_i,v_i\in \{0,1\}\}|.
$$
Notice that the distance may be zero even if the two strings are not the same. For instance, if $u=0\ast1\ast$ and $v=\ast 1\ast 0$, then $\dist(u,v)=0$. However, if $x,y\in \{0,1\}^n$, then $\dist(x,y)$ is equal to the Hamming distance between $x$ and $y$.

Clearly, $\dist(u,v)$ is exactly the number of dimensions in which the two corresponding boxes pass. Thus, by the above discussion, $n(k,d)$ is just the maximum number of strings in $S^d$ such that every two of them satisfy $1\leqslant \dist(u,v)\leqslant k$.

\begin{section}{Proof of Theorem \ref{main}}
Let us introduce some notation that will be helpful in our construction. For any two strings  $x=x_1x_2\cdots x_k$ and $y=y_1y_2\cdots y_l$,  the ordered pair $(x,y)$ is identified with their concatenation $xy=x_1x_2 \cdots x_ky_1y_2\cdots y_l$. Let us set $I=\{0,1\}$, $*^0=\{*\}$, $*^1=I$. If $s$ is a positive integer, then for any member  $\varepsilon$ of  the $s$-fold Cartesian product $I^s$  of $I$, we set $*^\varepsilon= *^{\varepsilon_1}\times *^{\varepsilon_2}\times\cdots\times *^{\varepsilon_s}$. If $x\in I^n$, then $\bar x\in I^n$ denotes the \textit{antipode} of $x$; that is $\bar x_i=1-x_i$, for $i=1,2,\ldots,n$. If $A\subseteq I^n$, then $\bar A=\{\bar x\colon x\in A\}$ is the antipode of the set  $A$. If $x, y\in I^n$, then $|x|=\sum_{i=1}^n x_i$ and $\langle x, y\rangle=\sum_i ^n x_iy_i$.

The \emph{diameter} $\diam(A)$ of $A\subseteq \{0,1,*\}^n$ is defined in the standard way, as the maximum distance $\dist(x,y)$ between any pair of elements $x,y\in A$. If $\diam A \leqslant k$, then we will call $A$ a \emph{$k$-neighborly code}. Let in addition $B\subseteq\{0,1,*\}^n$. It is convenient to consider the diameter of a pair $A$, $B$ of sets:
$$
\diam(A,B)=\max\{\dist(a,b)\colon a \in A, b\in B\}.
$$ 
Clearly, $\diam (A)=\diam(A,A)$.    

Let $\ka Z=\{Z_u\colon u \in U\}$ be an indexed family of subsets of a set $Z$. The family $\ka Z$ \textit{defines a partition} of $Z$ if $\bigcup \ka Z=Z$  and  for every pair $u, v\in U$, if $u\neq v$, then  $Z_u$ and $Z_v$ are disjoint. We do not rule out that  the empty set is a member of $\ka Z$.   

For $\varepsilon \in I^s$, if $\varepsilon \neq \alpha = (0,0,\ldots, 0)$ and $\varepsilon \neq \omega = (1,1,\ldots, 1)$, then there is a unique $i$ such that $\varepsilon_i\neq \varepsilon_{i+1}=\ldots =\varepsilon _s$. Let us set
$$
A_{\varepsilon}=\{\varepsilon_1\}\times\cdots\times\{\varepsilon_i\}\times I^{s-1-i}\subseteq I^{s-1}.
$$
Let us remark that the uniqueness of $i$ determines a function $t$ on the set of admissible $\varepsilon$'s: $\varepsilon \mapsto t(\varepsilon)=i$.  

\begin{lemma}
\label{part A}
 Let $s>1$, and $k\in \{1,2,\ldots, s-1\}$. The indexed family $\{A_\varepsilon\colon |\varepsilon|=k\}$ defines a partition of $I^{s-1}$. 
\end{lemma}
\proof Let $\varepsilon$ and $\gamma$ be different members of $I^s$ and let $|\varepsilon|=|\gamma|=k$. Let us fix the minimum index $i$ for which $\varepsilon_i\neq \gamma_i$. Suppose  
that $t(\varepsilon)< i$. Then $\varepsilon_{i}=\varepsilon_{i+1}=\ldots =\varepsilon_s$. Suppose first that $\varepsilon_{i}=0$. Then, by the definition of $i$, one has 
$
k=|\varepsilon|=\sum_{n=1}^{i-1} \varepsilon_n
$
 and 
 $
\gamma_j= \varepsilon_j,
 $
 whenever $j\leq i-1$. Since $\gamma_i$ has to be $1$, one gets $|\gamma|>k$, which is impossible. Suppose now that $\varepsilon_i= 1$. Then  
$$ 
k= |\varepsilon|= \sum_{n=1}^{i-1}\varepsilon_n +  1 +\sum_{n=i+1}^s 1> \sum_{n=1}^{i-1}\varepsilon_n +  0 +\sum_{n=i+1}^s \gamma_n =|\gamma|,
$$
which again is impossible   Consequently, $t(\varepsilon)\geqslant i$  and, by the same reason, $t(\gamma)\geqslant i$.  
 By the definition of $\varepsilon\mapsto A_\varepsilon$, it follows that the $i$-th factor of $A_\varepsilon$ is $\{\varepsilon_i\}$, while the corresponding factor of $A_\gamma$ is $\{\gamma_i\}$. Therefore, $A_\varepsilon$ and $A_\gamma$ are disjoint. 
 
 It remains to show that  $\{A_\varepsilon\colon |\varepsilon|=k\}$  is a covering of $I^{s-1}$. Let $x\in I^{s-1}$.  If $|x|<k$, then let us replace subsequent zeros appearing in $x$ by ones, starting from the right end, until we reach the element $u$ with the norm $|u|=k-1$.  Set $\varepsilon=(u,1)$. If $|x|\geq k$, then let us replace subsequent ones appearing in $x$ by zeros, starting from the right end, until we reach the element $v$ with the norm $|v|=k$. Now, set $\varepsilon=(v,0)$. It is easily seen that in each case under consideration $x$ is in the resulting $A_\varepsilon$. \hfill $\square$

 Let $d$ and $s$ be as in Theorem \ref{main}. For $k=0,1,\ldots, s$ we define sets $X_k\subseteq I^{d-s}$ as follows:
\begin{eqnarray*}
X_0 &=& \{x\in I^{d-s}\colon |x|\leqslant \frac d 2-s\}, \\
X_k &=& \{x\in I^{d-s}\colon |x|=\frac d  2 -s+k \}, \text{whenever $0<k<s$},\\
X_s &=& \left\{x\in I^{d-s}\colon |x|\geqslant \frac d 2\right\}.
\end{eqnarray*}
Clearly, these sets define a partition of the discrete box $I^{d-s}$. 
 
 For $\varepsilon \in I^s\setminus \{\alpha, \omega\}$,  let us set 
$$
 B_\varepsilon=\{x=(x', x'')\in I^{d-2s+1}\times I^{s-1}\colon x\in X_{|\varepsilon|} , x''\in A_\varepsilon\}
 $$
 Clearly, 
 \begin{equation}
 \label{inters}
 B_\varepsilon=X_{|\varepsilon|}\cap (I^{d-2s+1}\times A_\varepsilon)
 \end{equation}
 
 As an immediate consequence of Lemma \ref{part A}, we obtain
 \begin{corollary}
 \label{part B}
 Let $k\in \{1,2,\ldots, s-1\}$. The indexed family $\{B_\varepsilon\colon |\varepsilon|=k\}$ defines a partition of $X_k$. 
\end{corollary}

\begin{lemma}
\label{nier1}
Let $\varepsilon$ and $\gamma$ belong to $I^s\setminus\{\alpha, \omega\}$. Let $i=t(\varepsilon)$ and $j=t(\gamma)$. Let $\varepsilon_s=0$. If  one of the two possibilities
\begin{enumerate}
\item $\varepsilon_s=\gamma_s=0$,
\item $ \gamma_s=1$ and $j\geqslant i$ 
\end{enumerate}
holds true, then for every $x \in B_\varepsilon$ and $y\in B_\gamma$
$$
\skalar \varepsilon \gamma \leqslant \skalar x y.
$$

\end{lemma}
\proof  
 For a string $x=x_1 \cdots x_n$ and an index $k\leqslant n$, let $x|k=x_1\cdots x_k$. Then, by the assumptions on $\varepsilon$ and $\gamma$, we have
$$ 
\skalar \varepsilon \gamma=\skalar {\varepsilon|i} {\gamma|i}.
$$
We also have
$$
\skalar x y= \skalar {x'}{ y'} +\skalar{x''}{y''}\geqslant \skalar {x''|i}{y''|i},
$$
where $x=(x',x'')$, $y=(y',y'')$ are the decompositions of $x$ and $y$ as described in the definition of $\varepsilon \mapsto B_\varepsilon$. Since, in addition, $ \varepsilon|i \leqslant x''|i $ and $\gamma|i \leqslant y''|i $ with respect to the coordinatewise order, our inequality can be easily deduced. \hfill $\square$

For the elements $\alpha$ and $\omega$ of $I^s$, let us set $B_\alpha= X_{|\alpha|}=X_0$ and $B_\omega= X_{|\omega|}=X_{s}$. It is clear that each of the sets $X_0$ and $X_1$ is antipodal to the other. Therefore, 
$$
 \overline{B_\alpha}= B_\omega=B_{\bar \alpha}, \quad \text{and}  \quad \overline{B_\omega}=B_{\alpha}=B_{\bar\omega}.
$$
This observation is a part  of the following lemma:
\begin{lemma}
\label{antypody}
 For every positive integer $s$ and every $\varepsilon \in I^s$
 $$
 \overline{B_\varepsilon}=B_{\bar\varepsilon}.
 $$ 
\end{lemma}
\proof  Due to the comment preceding our lemma, we may assume that $|\varepsilon| \in \{1,\ldots, s-1\}$. It is easily seen that 
$\overline{X_{|\varepsilon|}}=X_{s-|\varepsilon|}= X_{|\bar\varepsilon|}$. 
Observe now that if $i=t(\varepsilon)$, then $i=t(\bar\varepsilon)$. Therefore, 
$$
A_\varepsilon=\{\varepsilon_1\}\times\cdots\times\{\varepsilon_i\}\times I^{s-1-i}\quad \text{and}\quad A_{\bar\varepsilon}=\{1-\varepsilon_1\}\times\cdots\times\{1-\varepsilon_i\}\times I^{s-1-i},
$$
which readily implies $\overline {A_\varepsilon}=A_{\bar\varepsilon}$. Consequently,
$$
\overline{B_\varepsilon}=\overline{X_{|\varepsilon|}\cap(I^{d-2s+1}\times A_\varepsilon})= \overline{X_{|\varepsilon|}}\cap (I^{d-2s+1}\times\overline{A_{\varepsilon}})= X_{|\bar\varepsilon|}\cap (I^{d-2s+1}\times A_{\bar\varepsilon} )= B_{\bar\varepsilon}.
$$
\hfill $\square$

\begin{lemma}
\label{diam}
For every $\varepsilon, \gamma \in I^s\setminus\{\alpha, \omega\}$
$$
\diam(B_\varepsilon, B_\gamma)\leqslant \min\{d-s-\skalar \varepsilon \gamma, d-s-\skalar {\bar\varepsilon} {\bar\gamma}\}.
$$
\end{lemma}

\proof 
Let $x \in B_\varepsilon$ and $y\in B_\gamma$. Since $x\in X_{|\varepsilon|}$ and $y\in X_{|\gamma|}$, we have
\begin{eqnarray*}
\dist(x,y) &=& |x|+|y|-2\skalar x y\\
&=& \frac d 2-s+|\varepsilon|+\frac d 2 -s+|\gamma|-2\skalar x y\\
&=& d-s +(|\varepsilon|+|\gamma|-s)-2\skalar x y.\\
\end{eqnarray*}
Since the expression in parentheses has a value not exceeding $\skalar \varepsilon \gamma$, we obtain
\begin{eqnarray*}
\dist(x,y) &\le& d-s +\skalar \varepsilon  \gamma - 2\skalar x y\\
&=&  d-s -\skalar \varepsilon  \gamma  +2(\skalar \varepsilon \gamma -\skalar x y).\\
\end{eqnarray*}
Suppose from now on that $\varepsilon$ and $\gamma$ are as assumed in Lemma \ref{nier1}. Then the expression $\skalar \varepsilon \gamma -\skalar x y)$ has non-positive value and
\begin{equation}
\label{one}
\dist(x,y)\leqslant d-s-\skalar \varepsilon \gamma.
\end{equation}
Similarly,
$$
\dist(x,y)=d-s -(s-|\varepsilon|-|\gamma|+\skalar \varepsilon \gamma) +(\skalar \varepsilon \gamma -2\skalar x y)
$$
The first term in parentheses is equal to $\skalar {\bar \varepsilon}{\bar\gamma}$, while the second one is, according to Lemma \ref{nier1}, non-positive, consequently
$$
\dist(x,y)\leqslant d-s- \skalar {\bar \varepsilon}{\bar\gamma},
$$
Combining (\ref{one}) with the last inequality, we  arrive to  the desired result. However, to complete the proof, we need to remove the constraints on $\varepsilon$ and $\gamma$. 

Note first that the constraints imposed in Lemma \ref{nier1} on $\varepsilon$ and $\gamma$ are not symmetric. If the pair $\varepsilon, \gamma$ does not satisfy them, then one of the pairs $\gamma$, $\varepsilon$ or $\bar\varepsilon, \bar\gamma$ or $\bar\gamma, \bar\varepsilon$ has to. It is obvious that it is enough to consider the case when the pair $\bar\varepsilon$, $\bar\gamma$ satisfies the assumptions of our lemma.  By Lemma \ref{diam}, the preceding part and the fact that $\varepsilon \mapsto \bar\varepsilon$ is an involution, we obtain
\begin{eqnarray*}
\diam(B_\varepsilon, B_\gamma) &= &\diam(\overline{B_\varepsilon}, \overline{B_\gamma}) \\
 &=& \diam(B_{\bar\varepsilon}, B_{\bar\gamma})\\
&=& \min\{d-s-\skalar {\bar\varepsilon} {\bar\gamma}, d-s-\skalar \varepsilon \gamma,\}. 
\end{eqnarray*}
\hfill $\square$ 

For $\varepsilon \in I^s$, let us set $C_\varepsilon=B_\varepsilon\times *^{\varepsilon}$. Since the members of the indexed family $\{B_\varepsilon\colon \varepsilon \in I^s\}$ are pairwise disjoint subsets of the discrete box $I^{d-s}$, the members of the family $\ka C=\{C_\varepsilon\colon \varepsilon \in I^s\}$ are also pairwise disjoint.  Moreover, $\bigcup \ka C$ is a $1$-neighborly code. In fact we have the following lemma.
 
\begin{lemma}
Let $d$ be a positive even integer and Let $s$ be a positive integer so that $2s\leqslant d$ The union of $\ka C =\{C_\varepsilon\colon \varepsilon \in I^s\}$ is a $(d-s)$-neighborly code. 
\end{lemma}
\proof 
Since $\ka C$ is a $1$-neighborly code, it suffices to show that for every $\varepsilon, \gamma \in I^s$, we have  
$ \diam(C_\varepsilon, C_\gamma)\leqslant d-s$. By the definition of $\ka C$,
$$
\diam(C_\varepsilon, C_\gamma)= \diam(B_\varepsilon, B_\gamma)+ \diam(*^\varepsilon, *^\gamma)
$$
As $\diam(*^\varepsilon, *^\gamma))=\skalar \varepsilon \gamma$, it follows from Lemma \ref{diam} that indeed  $\diam(C_\varepsilon, C_\gamma)$ does not exceed $d-s$ as long as $\varepsilon$ and $\gamma$ are different from $\alpha$ and $\omega$. 

 Suppose first that at least one of the two $\varepsilon$ and $\gamma$ is $\alpha$. Then $\diam(*^\varepsilon, *^\gamma))=0$. Moreover, $ \diam(B_\varepsilon, B_\gamma)\leqslant d-s$, as $B_\varepsilon$ and  $B_\gamma$ are subsets of $I^{d-s}$. Therefore,  also in this case $ \diam(C_\varepsilon, C_\gamma)\leqslant d-s$.

 Finally, suppose that one of the two $\varepsilon$ and $\gamma$ is $\omega$, while the other is in $I^{d-s}\setminus\{\alpha, \omega\}$. We may assume that $\varepsilon=\omega$. In this case,
 \begin{equation}
 \label{omega}
 \diam(*^\varepsilon, *^\gamma)= |\gamma|.
 \end{equation} 
 In order to estimate $D=\diam(B_\varepsilon, B_\gamma)$, we employ Lemma \ref{antypody}:
 $$
 D=\diam(B_\alpha, B_{\bar\gamma}).
 $$
 Let $x\in B_\alpha$ and $y\in B_{\bar\gamma}$.Since $B_\alpha \subseteq X_0$ and $B_{\bar\gamma}\subseteq X_{|\bar\gamma|}$, we have
$$
\dist(x,y)=|x-y|\leqslant |x|+|y|\leqslant \frac d 2 -s + \frac d 2-s+|\bar\gamma|=d-2s+(s-|\gamma|)= d-s-|\gamma| 
$$
 Consequently, $ D\leqslant d-s-|\gamma|$. This, combined with (\ref{omega}), implies the expected inequality. \hfill $\square$

\medskip
\noindent
\textit{Proof of Theorem \ref{main}.} Clearly, $n( d-s,d)\geqslant |\ka C|$. We have
$$
|\ka C|= \sum_\varepsilon |C_\varepsilon|= \sum_\varepsilon |B_\varepsilon|\cdot |*^{\varepsilon}|.
$$
Since $|*^{\varepsilon}|=2^{|\varepsilon|}$, we obtain
$$
|\ka C|= |B_\alpha| + \left(\sum_{k=1}^{s-1} 2^k\sum_{\varepsilon\colon |\varepsilon|=k} |B_\varepsilon|\right) +2^s|B_\omega|
$$
Since
\begin{itemize}
\item $B_\alpha= X_0$ and $|X_0|=\sum_{i=0}^{\frac d 2 -s-1} \binom {d-s} i$,
\item  $(B_\varepsilon\colon |\varepsilon|=k)$ is a partition of $X_k$, whenever $k=1,2,\ldots, s-1$ and $|X_k|=\binom {d-s}{\frac d 2 -s +k}$,
\item $B_\omega=X_s$ and $|X_s|=\sum_{i=\frac d 2 }^{d-s} \binom{d-s} i$, 
\end{itemize}
we easily come to the conclusion that the right side of the inequality under consideration expresses the cardinality of $\ka C$.\hfill$\square$ 

\end{section}

\begin{section}{Proof of Theorem \ref{asymptotics}}

	We shall assume for simplicity that both $d$ and $s$ are even. The remaining cases can be proven in a similar  manner.
	To get the left inequality, we can use our Theorem \ref{main}. In fact, we need less than that.  It is not so easy to prove Theorem 1 but it is almost trivial to show that  $C_\alpha\cup C_\omega$ is a $(d-s)$-neighborly. Thus, the cardinality of this set does not exceed $n(d-s, d)$. Consequently, 
	$$
	n(d-s, d)\geqslant   \sum_{i=0}^{\frac d 2 -s} \binom {d-s} i + 2^s\sum_{i=\frac d 2}^{d-s} \binom {d-s} i 
	$$
	By the symmetry of the binomial coefficients,
	$$ 
	n(d-s,d)\geqslant \frac 1 2 (1+2^s)\left (\sum_{i=0}^{\frac d 2 -s} \binom {d-s} i + \sum_{i=\frac d 2}^{d-s} \binom {d-s} i \right).
	$$
	Therefore,
	$$
	n(d-s,d)\geqslant \left(\frac 1 2 + 2^{s-1}\right)\left(2^{d-s} - \sum_{k=1}^{s-1} \binom {d-s}{\frac d 2 -s+k}\right).
	$$
	Since the central binomial coefficient is the largest and satisfies the inequality
	$$
	\binom{d-s}{\frac{d-s} 2} \leqslant \frac  {2^{d-s}}{\sqrt{\pi (d-s)/2}}, 
	$$
	we obtain
	$$
	n(d-s,d)\geqslant \left(\frac 1 2 +2^{s-1}\right)\left(2^{d-s} -  \frac  {2^{d-s}s}{\sqrt{\pi (d-s)/2}}\right). 
	$$
	Finally, 
	$$
	n(d-s,d)\geqslant \left( \frac 1 2 +\frac 1{2^{s+1}}\right)\left(1-\frac s {\sqrt{\pi (d-s)/2}}\right)2^d.
	$$
	
	The estimate on the right follows easily from the result by Cheng, Wang, Xu, and Yip [arxiv:2301.06485v1, Corollary 1.6], which we formulate below.
	\begin{theorem}
		Let $1\leqslant k\le d-2$ be integers. If $k=2t$, then 
		$$
		n(k,d)\leqslant 2^{d-2}+ 2^k+ \left(\frac 1 2 - \frac 1 {2^{d-k}}\right)\sum_{i=0}^t\binom d i.
		$$
		and if $k=2t+1$, then 
		$$
		n(k,d)\leqslant 2^{d-2}+ 2^k +\left( 1 -\frac 1 {2^{d-k-1}}\right)\sum_{i=0}^t \binom{d-1} i.
		$$
	\end{theorem}
	Observe first that $\sum_{i=0}^{t}\binom d i \leqslant \frac 1 2 2^d$ in the first case, while $\sum_{i=0}^{t} \binom{d-1} i \leqslant \frac 1 2  2^{d-1}$. Therefore, in both cases we have  
	$$
	n(k,d)\leqslant \frac 1 2 \left( 2^d + 2^k\right)
	$$
	and it remains to replace $k$ by $d-s$.
	
	Let us remark that the published version \cite{ChengWXY} of the paper by Cheng, Wang, Xu, and Yip (as well as its later arxiv versions) appears to be slightly different than the one we used. In particular, Corollary 1.6 from arxiv:2301.06485v1 is not reproduced there and it is not quite clear how it can be deduced without additional work.
\end{section}


\begin{thebibliography}{99}
	
	\bibitem{Aigner Ziegler}M. Aigner, G. Ziegler, Proofs from the Book, Sixth Edition, Springer Verlag, Berlin-Heidelberg, 2018.
	
	\bibitem{Alon} N. Alon,  Neighborly families of boxes and bipartite coverings, In: R. L. Graham at al. (eds.), \textit{The Mathematics of Paul Erd\"os II}, pp 27--31, Springer-Verlag, Berlin Heidelberg, 1997.
	
	\bibitem{Alon 2} N. Alon, Decomposition of the complete $r$-graph into complete $r$-partite $r$-graphs, \emph{Graphs Combin.} \textbf{2} (1986) 95--100.
	
	\bibitem{AlonGKP} N. Alon, J. Grytczuk, A.P. Kisielewicz and K. Przesławski, New bounds on the maximum number of neighborly boxes in $\mathbb{R}^d$, \textit{European J. Combin.} \textbf{114} (2023) 103797.
	
	\bibitem{ChengWXY} X. Cheng, M. Wang, Z. Xu, C. H. Yip, Exact values and improved bounds on $k$-neighborly families of boxes, \emph{European J. Combin.} $\mathbf{118}$ (2024) 103926.
	
	\bibitem{Cioaba} S. M. Cioabă, A. K\"{u}ngden, J. Verstra\"{a}te, On decompositions of complete hypergraphs, \emph{J. Combin. Theory Ser. A} \textbf{116} (2009) 1232--1234.
	
	\bibitem{GP} R. L. Graham and H. O. Pollak, On embedding graphs in squashed cubes, In: \textit{Lecture Notes in Mathematics} \textbf{303}, pp 99--110, Springer Verlag, New York-Berlin-Heidelberg, 1972.
	
	\bibitem{GKP} J. Grytczuk, A.P. Kisielewicz, K. Przes\l awski, Neighborly boxes and bipartite coverings; constructions and conjectures, arXiv:2402.02199.
	
	\bibitem{HS}  H. Huang and B. Sudakov, A counterexample to Alon-Saks-Seymour conjecture and related problems, \textit{Combinatorica} \textbf{32} (2012) 205--219.
	
	\bibitem{Luba} S. Łuba, Systems of Unit Cubes in $\mathbb{R}^d$, Master thesis, Uniwersytet Zielonogórski, 2023, (in Polish).
	
	\bibitem{Peck} G. W. Peck, A new proof of a theorem of Graham and Pollak, \textit{Discrete Math.} \textbf{49} (1984) 327--328.
	
	\bibitem{Tverberg} H. Tverberg, On the decomposition of $K_n$ into complete bipartite graphs, \textit{J. Graph Theory} \textbf{6} (1982) 493--494.
	
	\bibitem{van Lint} J.H. van Lint, $\{0,1,\ast\}$ distance problems in combinatorics, in: Surveys in Combinatorics 1985 (Glasgow, 1985), in: London Math. Soc. Lecture Note Ser., vol. 103, Cambridge Univ. Press, Cambridge, 1985, pp. 113--135.
	
	\bibitem{Vishwanathan1} S. Vishwanathan, A polynomial space proof of the Graham–Pollak theorem, \emph{J. Combin. Theory Ser. A} \textbf{115} (2008) 674--676.
	
	\bibitem{Vishwnathan2} S. Vishwanathan, A counting proof of the Graham–Pollak Theorem, \emph{Discrete Math.} \textbf{313} (2013) 765--766.
	
	\bibitem{Zaks3} J. Zaks, How Does a Complete Graph Split into Bipartite Graphs and How Are Neighborly Cubes Arranged?, \textit{Amer. Math. Monthly}  \textbf{92} (1985) 568--571.
	
\end{thebibliography}
\end{document}